\documentclass[a4paper,twoside,12pt]{article}
\usepackage[english]{babel}
\usepackage{bbding}
\usepackage[latin1]{inputenc}
\usepackage{kpfonts}
\usepackage{fancyhdr}
\usepackage{amsmath,amssymb,color,epsfig}
\usepackage{mathrsfs}
\usepackage{eufrak}
\usepackage{dsfont}
\usepackage{upgreek}
\usepackage{fancyhdr}
\usepackage[dvipsnames]{xcolor}
\newtheorem{theorem}{Theorem}[section]
\newtheorem{proposition}[theorem]{Proposition}

\def\1{\mathds{1}}
\def\e{\varepsilon}

\def\E{\mathbb{E}}

\title{Asymptotic density of zeros of half range generalized Hermite polynomials}

\date{}
\begin{document}
\author{Mohamed BOUALI}
\maketitle
\date
\begin{abstract}
We investigate the global density of zeros of generalized Hermite orthogonal polynomials,
subject to certain truncated conditions on its weight. We shall given explicitly the global density of zeros under some asymptotic conditions on the weight. Moreover we compute the asymptotic of the total energy of the equilibrium position of the system of $n$ movable
unit charges in an external field determined by the weight of the generalized Hermite polynomials. We will see that for finite $n$ the energy is in direct relationship with the zeros of the orthogonal polynomials.\\
{\bf Math Subject classification:} 15B52, 15B57, 60B10.
\end{abstract}
\section{Introduction}
Stieltjes \cite{st}, \cite{stl} considered the following electrostatic model. Fix two
charges $(\alpha+ 1)/2$ and $(\beta + 1)/2$ at $x = 1$ and $x = -1$, respectively, then
put $n$ movable unit charges at distinct points in $]-1; 1[$. The question is
to determine the equilibrium position of the movable charges when the interaction
forces arise from a logarithmic potential. Stieltjes proved that
the equilibrium position is attained at the zeros of the Jacobi polynomial
$P^{(\alpha;\beta)}_n(x)$. For a proof see {Szeg}o's book \cite{sz}. Another electrostatic problem
is to have a fixed point charge $(\alpha+1)/2$ at $x = 0$ and $n$ movable unit point
charges at distinct points in $[0;1[$. The state of equilibrium in the presence
of an additional external potential $v(x) = x$ is now reached at the zeros of
the Laguerre polynomial $L_n^{(\alpha)}(x)$ provided that the point charges interact
according to a logarithmic potential.

Let $a\geq 0$, $a<b$, and consider the electrostatic model, where we put an $n$ movable unit point charge at distinct points in $\Sigma_a=[a,+\infty[$ or $]-\infty,-a]\cup[a,+\infty[$. The state of the equilibrium in the presence of an additional external potential $V_n(x)=x^2+2\mu_n\log\frac1{|x|}+\log A_n(x)$ is realized at the zeros of the generalized Hermite polynomials $H_n^{\mu_n}(x, a)$ provided that the point charges interact
according to a logarithmic potential,
 where $H_n^{\mu_n}(x, a)$ are the orthogonal polynomials with respect the weight function $w_{\mu_n}(x)=C|x|^{2\lambda_n}e^{-x^2}$ restricted to $\Sigma_a$, $A_n(x)$ is some functions which will be given later and $C$ is a normalizing constant, see for instance \cite{I}. In other words, for $(x_1,...,x_n)\in\Sigma_a$, let $$E(x_1,...,x_n)=\sum_{i=1}^nV_n(x)+2\sum_{i<j}\log\frac1{|x_i-x_j|},$$ be the energy created by the $n$ movable unit charge $(x_1,...,x_n)$ in $\Sigma_a$. Then
 as it has been proved in \cite{I}, the equilibrium position is realized at $x_1^{(n, a)},...,x_n^{(n, a)}$ the zeros of $H_n^{\lambda_n}(x, a)$. Moreover
 $$E_n^*(a):=\min_{(x_1,...,x_n)\in \Sigma_a^n}E(x_1,...,x_n)=E(x_1^{(n, a)},...,x_n^{(n, a)}).$$
 One of the interesting questions is to study the asymptotic of $E_n^*(a)$ as $n$ goes to infinity. To do this, we consider on $\Sigma_a$ the probability measure
 $$\mu_n^a=\frac1n\sum_{i=1}^n\delta_{x_i^{(n, a)}}.$$
Using the general theory of logarithmic potential, we will prove that if $\lim_{n\to\infty}\frac{\lambda_n}{n}=\alpha\geq 0$, the measure $\mu_n^a$ scaled by the factor $\frac1{\sqrt n}$, converge tightly to some probability measure $\mu^a_\alpha$ with support  $S_{a, \alpha}\subset \Sigma_a$, which satisfies
$$\int_{S_{a, \alpha}}\log\frac1{|x-y|}\mu^a_\alpha(dy)+Q_\alpha(x)=\left\{\begin{aligned}&C\;{\rm if}\;x\in S_{a, \alpha}\\
&\geq C\;{\rm if}\;x\notin S_{a, \alpha}\end{aligned}\right.,$$
where $Q_\alpha(x)=x^2+2\alpha\log\frac1{|x|}$ and we give a explicit density to the measure above. Moreover we show that the asymptotic of the energy $E_n^*(a)$ is closely related to
$$E^*(a):=\int_{S_{a, \alpha}}\int_{S_{a, \alpha}}\log\frac1{|x-y|}\mu^a_\alpha(dy)\mu^a_\alpha(dx)+\int_{S_{a, \alpha}}Q_\alpha(x)\mu^a_\alpha(dx
).$$
The measure $\mu^a_\alpha$ is called the equilibrium measure, and $E^*(a)$ is the equilibrium energy. These results are closely related to those proved in \cite{bo} and \cite{boo}.

In the first section we define the system of orthogonal polynomials  in semi-finite interval of the form $[a,+\infty[$, and we give some simple asymptotic results for the coefficients of the the three terms recurrences relation. Moreover we give the second order differential equation satisfied by the  truncated generalized Hermite polynomials which has been studied in \cite{IC}, \cite{I}. In section 2 we defined the electrostatic energy and we show that the minimum of the energy is attained at the zeros of the generalized Hermite polynomials. Furthermore we give the central result, which state the convergence of the global density of zeros to some probability density. This is the subject of theorem 2.3. We end the section by an asymptotic formula of the energy, where we prove that $\displaystyle E_n^*(a)=\E^*n^2-n\Big(\lambda_n+\frac n2\Big)\log n+o(n^2),$ where $\E^*$ is some constant and $o(x)$ is a small terms in $x$ . For $\lambda_n=0$,   $\displaystyle E_n^*(a)\sim-\frac {n^2}2\log n,$ where $a_n\sim b_n$ as $n\to +\infty$  means $a_n/b_n\to 1$ as $n\to +\infty$. Such result has been proved in \cite{I}, and references therein. Section 4 is the object to defined the truncated generalized Hermite polynomials where we restrict the weight to the symmetric set $]-\infty,a]\cup[a,+\infty[$ $(a\geq 0)$, moreover we give some asymptotic of the coefficient of the three terms recurrences relation, which will be used later.  As in the previous section we will give the global density of zeros when $n$ goes to infinity, this is the object of theorem 5.2. In the last of section 3 and 4 we give some plots, which completes our analysis about the global density of zeros.
 \section{Half range generalized Hermite polynomials}
Let $a\in\Bbb R$ and $\lambda> -\frac12$.
Consider on $[a,+\infty[$ the normalizing weight $$\displaystyle w_\lambda(x)=\frac{1}{C_\lambda}|x|^{2\lambda}e^{-x^2},$$ where  $\displaystyle C_\lambda=\int_a^{+\infty}|x|^{2\lambda}e^{-x^2}dx$. For $a=0$, $\displaystyle C_\lambda=\frac{\Gamma(\lambda+\frac12)}{2}.$

On the space $\cal P$ of polynomials in one variable with real coefficients we
consider the inner product
$$\langle p,q\rangle=\int_a^{+\infty}p(x)q(x)w_\lambda(x)dx.$$
which makes $\cal P$ into a preHilbert space. From the system $(t^m, m \in\Bbb N)$
the Schmidt orthonormalizing produces a sequence $(H^{\lambda, a}_m)_m$ of orthogonal
polynomials: $H^{\lambda, a}_m$ is a polynomial of degree $m$ and
$$\int_a^{+\infty}H^{\lambda, a}_m(x)H^{\lambda, a}_n(x)w_\lambda(x)dx=\delta_{nm}.$$
$H^{\lambda, a}_m$ is said the half range generalized Hermite polynomials.

We know that the sequence $H^{\lambda, a}_n$ satisfies three-terms recurrence relation $\displaystyle H_0^{\lambda, a}(x)=1$, $\displaystyle H_1^{\lambda, a}(x)=\frac{(x-a_0)}{b_1}$,
$$xH^{\lambda, a}_n(x)=a_{n+1}H^{\lambda, a}_{n+1}(x)+b_nH^{\lambda, a}_n(x)+a_nH^{\lambda, a}_{n-1}(x),$$
where $b_n=\langle xH^{\lambda, a}_n,H^{\lambda, a}_n\rangle$ and $a_n=\langle xH^{\lambda, a}_n,H^{\lambda, a}_n\rangle$. Moreover if we write $H^{\lambda, a}_n(x)=\gamma_nx^n+$ lower order terms, then $\displaystyle a_n=\frac{\gamma_{n-1}}{\gamma_n}$.

The sequences $a_n, b_n$  and $\gamma_n$ all depend in the parameters $\lambda$ and $a$, but we omit this dependance in the notations.

\begin{proposition}---
for all $n\in\Bbb N$,
\begin{description}
\item[(1)] $$b_n\geq \frac12(H^{\lambda, a}_n(a))^2w_\lambda(a),$$
\item[(2)]$$a_{n+1}^2+a_n^2+b_n^2=n+\lambda+\frac12+\frac12a(H^{\lambda, a}_n(a))^2w_\lambda(a),$$
\end{description}
moreover for $n$ large enough and for all $\lambda\geq 0$, $b_n=O\big(\sqrt {n+\lambda}\big)$ and, $(H^{\lambda, a}_n(a))^2w_\lambda(a)=O\big(\sqrt{n+\lambda}\big)$.
\end{proposition}
{\bf Proof.}---

{\bf Step (1):} By orthogonality we know
$$\int_a^{+\infty}H^{\lambda, a}_n(x)'H^{\lambda, a}_n(x)w_\lambda(x)dx=0,$$
furthermore, integrate by part
$$\int_a^{+\infty}H^{\lambda, a}_n(x)'H^{\lambda, a}_n(x)w_\lambda(x)dx=-\frac12(H^{\lambda, a}_n(a))^2w_\lambda(a)+b_n
-\frac{\lambda}{C_\lambda}\int_a^{+\infty}(H^{\lambda, a}_n(x))^2x^{2\lambda-1}e^{-x^2}dx,$$
it follows that $$ b_n\geq\frac12(H^{\lambda, a}_n(a))^2w_\lambda(a).$$
{\bf Step (2):}
Integrate by part it follows that
$$\begin{aligned}&\int_a^{+\infty}H^{\lambda, a}_n(x)'H^{\lambda, a}_n(x)xw_\lambda(x)dx
=\langle xH^{\lambda, a}_n,xH^{\lambda, a}_n\rangle-(\lambda+\frac12)||H^{\lambda, a}_n||^2-\frac12a(H^{\lambda, a}_n(a))^2w_\lambda(a),\end{aligned}$$
moreover $$xH^{\lambda, a}_n(x)'=nH^{\lambda, a}_n(x)+c_{n-1}H^{\lambda, a}_{n-1}(x)+...+c_0H_0^\lambda(x),$$
since $$\int_a^{+\infty}H^{\lambda, a}_n(x)'H^{\lambda, a}_n(x)xw_\lambda(x)dx=\langle x(H^{\lambda, a}_n)',H^{\lambda, a}_n\rangle=n||H^{\lambda, a}_n||^2=n.$$
And  from the recurrence relation $$\langle xH^{\lambda, a}_n,xH^{\lambda, a}_n\rangle=a_{n+1}^2+a_n^2+b_n^2,$$
hence $$a_{n+1}^2+a_n^2+b_n^2=n+\lambda+\frac12+\frac12a(H^{\lambda, a}_n(a))^2w_\lambda(a).$$
To see the boundedness of $b_n$, using step (1) and (2), it follows that
$$b_n^2\leq n+\lambda+\frac12+|a|b_n,$$
hence $$(b_n-\frac12|a|)^2\leq  n+\lambda+\frac12+\frac14a^2,$$
which means that $$0\leq b_n\leq\sqrt{n+\lambda+\frac12+\frac14a^2}+\frac{|a|}{2}.$$
And $\displaystyle\frac{b_n}{\sqrt{n+\lambda}}$ is bounded for all $n$.
This complete the proof.\qquad\qquad\qquad\qquad\qquad\qquad\qquad$\blacksquare$ \\
\begin{proposition}---
The polynomials $H^{\lambda, a}_n$ satisfies the following differential equations
\begin{description}
\item[(1)]
$$H^{\lambda, a}_n(x)'=A_n(x)H^{\lambda, a}_{n-1}(x)-B_n(x)H^{\lambda, a}_n(x),$$
where $$A_n(x)=2a_n+\frac{2a_nb_n}{x}+\frac{a_nw_\lambda(a)(H^{\lambda, a}_n(a))^2}{x-a},$$ and $$ B_n(x)=\frac{a_nw_\lambda(a)H^{\lambda, a}_n(a)H^{\lambda, a}_{n-1}(a)}{x-a}+\frac{2a_n^2-n-a_nw_\lambda(a)H^{\lambda, a}_n(a)H^{\lambda, a}_{n-1}(a)}{x}.$$
\item[(2)]
$$H^{\lambda, a}_n(x)''+R_n(x)H^{\lambda, a}_n(x)'+S_n(x)H^{\lambda, a}_n(x)=0,$$
where $$R_n(x)=-2x+\frac{2\lambda}{x}-\frac{A_n'(x)}{A_n(x)},$$ and
$$S_n(x)=B'_n(x)-B_n(x)\frac{A'_n(x)}{A_n(x)}-B_n(x)\Big(2x-\frac{2\lambda}{x}+B_n(x)\Big)+\frac{a_n}{a_{n-1}}A_n(x)A_{n-1}(x).$$
\end{description}
\end{proposition}
See \cite{B}, \cite{C} and \cite{IC} for the proof.
\section{ Density of zeros of half range Hermite polynomials}
We propose that the weight function $w(x)$ creates two external fields. One is a
long range field whose potential at a point $x>a$, is $$\displaystyle -\log w(x)=x^2+2\lambda\log\frac1x+\log C_\lambda.$$
In addition in
the presence of $n$ unit charges, $w$ produces a short range field whose potential
is $\displaystyle\log\Big(\frac{A_n(x)}{a_n}\Big)$. Thus the total external potential $V (x)$ is the sum of the
short and long range potentials, that is
$$V(x)=x^2+2\lambda\log\frac1x+\log C_\lambda+\log\Big(2+\frac{2b_n}{x}+\frac{w_\lambda(a)(H^{\lambda, a}_n(a))^2}{x-a}\Big).$$
Consider the system
of $n$ movable unit charges in $[a; +\infty]^n$. In the presence of the external potential $V(x)$. If $x_1,...,x_n$ are the position of the particles arranged in decreasing order. The total energy of the system is
 $$ F(x_1,\cdots,x_n)=2\sum_{i< j}\log\frac{1}{|x_i-x_j|}+\sum_{i=1}^nV(x_i).$$
 \begin{proposition}---The energy $F(x_1,...,x_n)$ is minimal at the zeros of half range Hermite polynomials $H^{\lambda, a}_n$.
 \end{proposition}
 See for instance \cite{I}.

 Let $x^{(n, a)}_1,...,x^{(n, a)}_n\in[a, +\infty[$ denote the zeros of half range generalized Hermite polynomials $H^{\lambda_n, a}_n$, where $\lambda_n$ be some positive real sequence, and defined on $[a, +\infty[$  the probability measure
$$\mu_n^a=\frac1n\sum_{i=1}^n\delta_{\sigma^{(n, a)}_i}, $$
where $\displaystyle \sigma^{(n, a)}_i=\frac{x^{(n, a)}_i}{\sqrt n}$ are the rescaled zeros. Moreover  consider the modified rescaled energy in the sense
 $$E(x_1,\cdots,x_n)=F(x_1,...,x_n)-n\log C_{\lambda_n}=2\sum_{i< j}\log\frac{1}{|x_i-x_j|}+n\sum_{i=1}^nV_n(x_i),$$
 where $\displaystyle V_n(x)=x^2+2\alpha_n\log\frac1x+\frac1n\log\Big(2+\frac{2b_n}{\sqrt n x}+\frac{\beta_n}{\sqrt nx-a}\Big)$, $\displaystyle\alpha_n=\frac{\lambda_n}{n}$, and $\displaystyle\beta_n=w_{\lambda_n}(a)(H_n^{\lambda_n}(a))^2$. Then the minimum of the energy is attained at the rescaled zeros $\displaystyle \sigma_i^{(n, a)}$ of the polynomials $H_n^{\lambda_n, a}$.

{\bf Through all the rest of the paper it will be assumed that $a>0$ if $\alpha\geq 0$, and $a\in\Bbb R$ if $\alpha=0$, where $\alpha=\lim\limits_{n\to +\infty}\alpha_n$.}
\begin{theorem}---
Assume that the sequence $\alpha_n$ converges to some $\alpha\geq 0$. Then the measure $\mu_n^a$ converge for the tight topology to some probability measure $\mu_\alpha^\sigma $ supported by $[\sigma,b]$ with density
$f_{\alpha, \sigma}$, where $\sigma=\sigma(a, \alpha)$, and such that:
 \begin{description}
\item[(1)] For $\alpha>0$, and $a=0$, there is a unique $\sigma_0=\sigma(\alpha)>0$, and a unique $b>\sigma_0$, and the density is $$f_{\alpha, \sigma_0}=\frac{1}{\pi}\sqrt{(b-x)(x-\sigma_0)}\Big(1+\frac{\alpha}{\sqrt{\sigma_0 b}}\frac1x\Big),$$
    where $$\sigma_0=\sqrt{\frac53+\frac{5\alpha}{3}-\frac\beta3-\frac23\sqrt{2+4\alpha-4\alpha^2+2(1+\alpha)\beta}},$$
 $\displaystyle\beta=\sqrt{1+2\alpha+4\alpha^2}.$ And
 $$\displaystyle b=b(\alpha)=\frac23\Big(\sqrt{6(\alpha+1)-2\sigma_0^2}-\frac{\sigma_0}2\Big).$$
\item[(2)] For $\alpha\geq 0$, and $a>0$, then $\sigma=a$, and
$$f_{\alpha, \sigma}=\frac{1}{2\pi}\sqrt{\frac{b-x}{x-\sigma}}\Big(2x+b-\sigma-2\alpha\sqrt{\frac{\sigma}{b}}\frac1x\Big),$$
where $\sigma<b$, and $b$ is the unique solution of the following equations
$$\frac{3}{4}(b-\sigma)^2+\sigma(b-\sigma)+2\alpha\sqrt{\frac \sigma b}-2\alpha-2=0,\;{\rm and\;} b+\sigma-\frac{2\alpha}{\sqrt{\sigma b}}\geq 0.$$
\item[(3)] For $\alpha=0$, and $a\in\Bbb R$, two cases are present:
\item[(a)] if $a\geq -\sqrt 2$, then $\sigma=a$.

\item [(b)] if $a<-\sqrt 2$ then $\sigma=-\sqrt 2$.

In each cases
$$f_{0, \sigma}=\frac{1}{2\pi}\sqrt{\frac{b-x}{x-\sigma}}\Big(2x+b-\sigma\Big),$$
and $$b=\frac23\Big(\sqrt{\sigma^2+6}+\frac \sigma2\Big).$$
\end{description}
\end{theorem}
From {\bf (3)} if we put $\alpha=a=0$, one obtains $$f_{0,0}(x)=\frac{1}{2\pi}\sqrt{\frac{b-x}{x}}\Big(2x+b\Big),$$
where $\displaystyle b=\frac23\sqrt{6}$. Such density appear in first time in \cite{S}, where the authors steadies the probability that all eigenvalues of Gaussian Hermitian matrix are  positives. A generalization of such question is treated in \cite{bo}. The author study the asymptotic density of eigenvalues in the interval $[a,+\infty[$ $(a\in\Bbb R)$ for the Generalized Gaussian unitary ensemble matrices. He show that the probability that all eigenvalues of Generalized Gaussian Hermitian matrix to be positives is given by $f_{\alpha, a}$.

In the third case on can see that for $a\leq 0$, $f_{0, a}$ is a probability density, if $b\geq -a$, hence
$\sqrt{a^2+6}\geq -2a$, which give that $a\geq-\sqrt 2$. Which means that the limit case for which $f_{\alpha, a}$ is a probability density is $a\geq -\sqrt 2$, and the limit density in that case is the semi-circle law
$$f_{0,-\sqrt 2}(x)=\frac{1}{\pi}\sqrt{2-x^2}.$$
In this case $f_{0,-\sqrt 2}$ is the distribution of zeros of the classical hermite polynomials $H^0_n$. Moreover one can prove by potential theory that for all $a<-\sqrt 2$, the density of zeros still unchanged, and is described by the semi-circle law. Which explain that the density of the zeros of $H^0_n$, when we restricted the weight $w_0$ to the half line $[a,+\infty[$, for $a\leq -\sqrt 2$, approximate the density of zeros of Hermite polynomials $H_n$, for the weight function $w_0$ on all the real line.

To prove the theorem we need some preliminary results.

Let ${\mathfrak M}^1(\Sigma)$ be the set of probability
measures on the closed set $\Sigma\subset\Bbb R$. We equip the space ${\mathfrak M}^1(\Sigma)$ with the tight topology. For this
topology a sequence $(\nu_n)$ converge to a measure $\nu$ if, for every continuous
bounded function $f$ on $\Sigma$,
$$\lim_{n\to\infty}\int_{\Sigma}f(x)\nu_n(dx)=\int_{\Sigma}f(x)\nu(dx).$$
This topology is metrizable. If $\Sigma$ is bounded, then ${\mathfrak M}^1(\Sigma)$ is compact.\\
Let $\Sigma$ be a closed set $\Big(\Sigma = \Bbb R,\; ]-\infty,a],\; [ a, +\infty [ \;{\rm or}\; ]-\infty,a]\cup[ a, +\infty [\Big)$, and $Q$ a
function defined on $\Sigma$ with values on $] -\infty,+\infty]$, continuous on int$(\Sigma)$. If
$\Sigma$ is unbounded, it is assumed that
$$\lim_{|x|\to +\infty}\Big(Q(x)-\log(1+x^2)\Big)=\infty.$$
If $\nu$ is a probability measure supported by $\Sigma$, the energy $E(\nu)$ of $\nu$ is
defined by
$$\E(\nu)=\int_{\Sigma}U^{\nu}(x)\nu(dx)+\int_{\Sigma}Q(x)\nu(dx),$$
where $$U^{\nu}(x)=\int_{\Sigma}\log\frac{1}{|x-t|}\nu(dt).$$
 By a straightforward computation one can prove that $E(\nu)$ is bounded from below. Hence we defined
 $$\E^*=\inf\Big\{\E(\nu)\mid \nu\in\mathfrak{M}^1(\Sigma)\Big\}.$$
 \begin{theorem}---
If $\nu(dx)=f(x)dx$, where $f$ is a continuous function with
compact support $\subset\Sigma$. Then the potential $\displaystyle U^{\nu}(x)=\int_\Sigma\log\frac{1}{|x-t|}\nu(dt)$ is a continuous function, and,
$\E^*\leq \E(\nu) < \infty$. Furthermore there is a unique measure $\nu^*\in\mathfrak{M}^1(\Sigma)$ such that
$$\E^*=\E(\nu^*).$$
The support of $\nu^*$ is compact.\\
This measure $\nu^*$ is called the equilibrium measure.
\end{theorem}
\begin{proposition}---
Let $\nu\in\mathfrak{ M}^1(\Sigma)$ with compact support.
Assume that the potential $U^{\nu}$ of $\nu$ is continuous and that there is a
constant $C$ such that\\
${\rm (i)}$ $U^{\nu}(x)+\frac 12Q(x)\geq C$ on $\Sigma$.\\
${\rm (ii)}$ $U^{\nu}(x)+\frac 12Q(x)= C$ on ${\rm supp}(\nu)$. Then $\nu$ is the equilibrium measure: $\nu=\nu^*$.
\end{proposition}
The constant $C$ is called the (modified) Robin constant. Observe that
$$\E^*=C+\frac 12\int_{\Sigma}Q(x)\nu^*(dx).$$
See for the proofs of the previous theorem and proposition \cite{F}.

The proof of the next proposition can be found in \cite{bo}.
\begin{proposition}---
 Consider the case where $\displaystyle\Sigma=[a,+\infty[$, $\displaystyle Q_{\alpha}(s)=s^2+2\alpha\log\frac{1}{s}$, with $\alpha\geq 0$. Then the equilibrium measure is $\mu_\alpha^a$ of the theorem {\rm(2.2)}. Moreover for $\alpha=0$, the energy is given by
$$\E^*= \frac{1}{108} \Big(81+72 a^2-2 a^4+\big(30a +2 a^3\big) \sqrt{6+a^2}-108 \log \Big(\frac{1}{6} \big(-a+\sqrt{6+a^2}\big)\Big)\Big).$$
 and for $\alpha=a=0$,
 $$\E^*=\frac34+\frac12\log6.$$
 \end{proposition}
 For the proof of the value of the energy in the case $\alpha=0$, see for instance \cite{S}.

{\bf Proof of theorem 2.2}---
We denote in the proof the equilibrium energy by $$\E^*_\alpha=\E(\mu_\alpha^a).$$ Let defined $$\tau_n=\frac{1}{n(n-1)}E^*_{\alpha_n, n}=\frac{1}{n(n-1)}E(\sigma_1^{(n, a)},...,\sigma_n^{(n, a)}).$$
Consider for a probability measure $\mu$  the energy
$$\E_{\alpha_n}(\mu)=\int_a^{+\infty}\int_a^{+\infty}\log\frac{1}{|s-t|}\mu(ds)\mu(dt)+\int_a^{+\infty }Q_{\alpha_n}(s)(s)\mu(ds),$$
where $\displaystyle Q_{\alpha_n}(s)=s^2+2\alpha_n\log\frac1s$.

Since for a probability $\mu$ with support in $[a,+\infty[$
$$\begin{aligned}&\int_{[a,+\infty[^n}E(x_1,...,x_n)\mu(dx_1)...\mu(dx_n)\\&=n(n-1)\int_a^{+\infty}\int_a^{+\infty}\log\frac{1}{|s-t|}\mu(ds)\mu(dt)+
n^2\int_a^{+\infty}Q_{\alpha_n}(s)\mu(ds)\\&+n\int_a^{+\infty}\log\Big(2+\frac{2b_n}{\sqrt n s}+\frac{\beta_n}{\sqrt n s-a}\Big)\mu(ds).\end{aligned}$$
Hence $$\begin{aligned}\tau_n&\leq\E_{\alpha}(\mu)+\frac{n^2}{n(n-1)}\int_a^{+\infty}\Big(Q_{\alpha_n}(s)-Q_\alpha(s)\Big)\mu(ds)
+\frac{1}{n-1}\int_a^{+\infty}Q_{\alpha}(s)\mu(ds)\\&+\frac1{n-1}\int_a^{+\infty}\log\Big(2+\frac{2b_n}{\sqrt n s}+\frac{\beta_n}{\sqrt n s-a}\Big)\mu(ds).\end{aligned}$$
For $\mu=\mu_\alpha^a$ the equilibrium measure, one gets
\begin{equation}\begin{aligned}\label{eq1}\tau_n&\leq\E^*_{\alpha}+\frac{n^2}{n(n-1)}\int_a^{+\infty}\Big(Q_{\alpha_n}(s)-Q_\alpha(s)\Big)
\mu_\alpha^a(ds)+\frac{1}{n-1}\int_a^{+\infty}Q_{\alpha}(s)\mu^a_\alpha(ds)\\&+\frac1{n-1}\int_a^{+\infty}\log\Big(2+\frac{2b_n}{\sqrt n s}+\frac{\beta_n}{\sqrt n s-a}\Big)\mu_\alpha^a(ds).\end{aligned}\end{equation}
Since $$\Big|\int_a^{+\infty}(Q_{\alpha_n}(s)-Q_{\alpha}(s))\mu_\alpha^a(ds)\Big|\leq 2C_0|\alpha_n-\alpha|,$$
where $\displaystyle C_0=\int_a^{+\infty}|\log(x)|\mu_\alpha^a(dx).$
Which implies that \begin{equation}\label{eq2}\lim_{n\to+\infty}\int_a^{+\infty}(Q_{\alpha_n}(s)-Q_c(s))\mu_\alpha^a(ds)=0.\end{equation}
From proposition (2.1) and the fact that $\alpha_n$ converge then
$\displaystyle \frac{b_n}{\sqrt {n}}$ and $\displaystyle \frac{\beta_n}{\sqrt {n}}$ are bounded by some positive constants $c_2$ and $c'_2$, and then
$$0\leq\frac1{n-1}\int_a^{+\infty}\log\Big(2+\frac{b_n}{\sqrt n s}+\frac{\beta_n}{\sqrt n s-a}\Big)\mu_\alpha^a(ds)\leq\frac1{n-1}\int_a^{+\infty}\log\Big(2+\frac{2c_2}{s}+\frac{c'_2}{ s-a}\Big)\mu_\alpha^a(ds),$$
since the integral $\displaystyle\int_a^{+\infty}\log\Big(2+\frac{2c_2}{s}+\frac{c'_2}{ s-a}\Big)\mu_\alpha^a(ds)$ converge, hence
\begin{equation}\label{eq4}\lim_{n\to+\infty}\frac1{n-1}\int_a^{+\infty}\log\Big(2+\frac{b_n}{\sqrt n s}+\frac{\beta_n}{\sqrt n s-a}\Big)\mu_{\alpha,a}(ds)=0.\end{equation}
From equations (\ref{eq1}), (\ref{eq2}) and (\ref{eq4}), one gets
\begin{equation}\label{eq5}\limsup_n\tau_n\leq\E^*_\alpha.\end{equation}
Let $$k_n(s,t)=\log\frac{1}{|s-t|}+\frac12V_n(s)+\frac12V_n(t),$$
and $$h_n(t)=V_n(t)-\log(1+t^2),$$
where  $\displaystyle V_n(x)=x^2+2\alpha_n\log\frac1x+\frac1n\log\Big(2+\frac{2b_n}{\sqrt n x}+\frac{\beta_n}{\sqrt n x-a}\Big)$.
The positive sequence $\alpha_n$ convergent, hence there is some positive constants $a_1, a_2\geq 0$, such that $a_1\leq\alpha_n\leq a_2$. From the fact that for all $x>a$,\\ $\displaystyle \frac1n\log\Big(2+\frac{2b_n}{\sqrt n x}+\frac{\beta_n}{\sqrt n x-a}\Big)\geq 0$ one gets,
$$h_n(t)\geq \left\{\begin{aligned}&t^2+2a_1\log\frac1t-\log(1+t^2)=h_1(t),\quad{\rm if}\;a<t\leq 1\\&t^2+2a_2\log\frac1t-\log(1+t^2)=h_2(t),\quad{\rm if}\;t\geq 1.\end{aligned}\right. $$
Let $$h(t)=\inf(h_1(t), h_2(t)),$$ and
using the inequality $$|s-t|\leq\sqrt{1+s^2}\sqrt{1+t^2},$$ it follows that, for all $n\in\Bbb N$ and all $s,t>0$,
$$k_n(s,t)\geq \frac12h(s)+\frac12h(t).$$
 and
$$\sum_{i\ne j}k_n(x_i,x_j)\geq (n-1)\sum_{i=1}^nh(x_i),$$
hence $$ \sum_{i\neq j}\log\frac{1}{|x_i-x_j|}+(n-1)\sum_{i=1}^nV_n(x_i)\geq (n-1)\sum_{i=1}^n h(x_i),$$
furthermore $$V_n(x)\geq h(x),$$
hence
$$E^*_n=E\Big(\sigma^{(n, a)}_1,...,\sigma^{(n, a)}_n\Big)\geq n^2\int_a^{+\infty}h(t)\mu_n^a(dt).$$
From (\ref{eq1}) one gets,
$$\int_a^{+\infty}h(t)\mu_n^a(dt)\leq\frac{n(n-1)}{n^2}\tau_n\leq \frac{n-1}{n}\big(\E^*_\alpha+B_n\big),$$
where $$\begin{aligned}B_n&=\frac{n^2}{n(n-1)}\int_a^{+\infty}\Big(Q_{\alpha_n}(s)-Q_\alpha(s)\Big)
\mu_\alpha^a(ds)+\frac{1}{n-1}\int_a^{+\infty}Q_{\alpha}(s)\mu^a_\alpha(ds)\\&+\frac1{n-1}\int_a^{+\infty}\log\Big(2+\frac{2b_n}{\sqrt n s}+\frac{\beta_n}{\sqrt n s-a}\Big)\mu_\alpha^a(ds).\end{aligned}$$
It has be seen in the previous that $B_n$ goes to $0$ as $n\to +\infty$, thus, there is some constant $C_0$ such that for all $n\in\Bbb N$,
$$\int_a^{+\infty}h(t)\mu_n^a(dt)\leq C_0.$$
Moreover $\lim\limits_{t\to+\infty}h(t)=+\infty$, then by the Prokhorov criterium there is some subsequence $n_j\to\infty$ such that the measure $\mu^a_{n_j}$ converge for the tight topology to some probability measure $\sigma^a$. We will denote simply $\mu_n^a$ the subsequence.

For $\ell>0$ consider the cut kernel $\displaystyle k_n^\ell(s,t)=\inf(k_n(s,t),\ell),$
and let $$\widetilde{k_{\alpha_n}}(s,t)=\log\frac{1}{|s-t|}+\frac12Q_{\alpha_n}(s)+\frac12Q_{\alpha_n}(t),$$
where $Q_{\alpha_n}(s)=s^2+2\alpha_n\log\frac1s,$ and
 $\displaystyle \widetilde{k_n}^\ell(s,t)=\inf(\widetilde{k_n}(s,t),\ell).$
It is easy to see that for all $n\in\Bbb N$, and all $s,t\in\Bbb R^*_+$,
$$\widetilde{k_{\alpha_n}}(s,t)\leq k_n(s,t),$$
and
\begin{equation}\label{eq6}\widetilde{k_{\alpha_n}}^\ell(s,t)\leq k_n^\ell(s,t).\end{equation}
Let $\e>0$, there is $n_0$, such that for all $n\geq n_0$,
$$\alpha-\e\leq\alpha_n\leq \alpha+\e,$$
Let $n\geq n_0$, dived $\Bbb R^2_+\setminus\{(s,t)\mid s= t\,{\rm and }\;s= 0,\,{\rm and}\; t=0\} $ to four region
$$R_1=\{(s,t)\mid s\geq 1\,{\rm and}\; t\ge 1\},\quad R_2=\{(s,t)\mid a<s\leq 1\,{\rm and}\; a<t\leq 1\},$$
and $$R_3=\{(s,t)\mid a<s\leq 1\,{\rm and}\; t\geq 1\},\quad R_4=\{(s,t)\mid \,s\geq 1\,{\rm and}\; a<t\leq 1\}.$$
If $(s,t)\in R_1$, then
$$\widetilde{k_{\alpha_n}}(s, t)\geq \widetilde{k_{\alpha+\e}}(s,t).$$
If $(s,t)\in R_2$,
$$\widetilde{k_{\alpha_n}}(s,t)\geq \widetilde{k_{\alpha-\e}}(s,t).$$
If $(s,t)\in R_3$,
$$\widetilde{k_{\alpha_n}}(s,t)\geq \log\frac{1}{|s-t|}+\frac12Q_{\alpha+\e}(t)+\frac 12Q_{\alpha-\e}(s),$$
hence $$\widetilde{k_{\alpha_n}}(s,t)\geq\frac12\left(\widetilde{k_{\alpha+\e}}(s,t)+\widetilde{k_{\alpha-\e}}(s,t)\right).$$
By symmetry of the kernel $\widetilde{k_{\alpha_n}}$  the last inequality  is valid in $R_4$.\\
We obtain in $\Bbb R^2_+\setminus\Big\{(s,t)\mid s= t\,{\rm and }\;s= 0,\,{\rm and}\; t=0\Big\} $,
$$\widetilde{k_{\alpha_n}}(s, t)\geq  \theta_1\,\widetilde{k_{\alpha+\e}}(s,t)+\theta_2\,\widetilde{k_{\alpha-\e}}(s,t),$$
where $(\theta_1,\theta_2)=(1,0)$ in $R_1$, $(\theta_1,\theta_2)=(0,1)$ in $R_2$ and $(\theta_1,\theta_2)=(\frac12,\frac12)$ in $R_3\cup R_4$.

Hence if we take the infimum then for all $(s,t)\in]a,+\infty[^2$, $$\widetilde{k_{\alpha_n}}^\ell(s, t)\geq \theta_1\widetilde{k_{\alpha+\e}}^\ell(s,t)+\theta_2\widetilde{k_{\alpha-\e}}^\ell(s,t),$$
By making use of equation(\ref{eq6}) it yields
$$k_n^\ell(s, t)\geq \theta_1\widetilde{k_{\alpha+\e}}^\ell(s,t)+\theta_2\widetilde{k_{\alpha-\e}}^\ell(s,t).$$
Take the energy in both sides of the previous inequality, it follows for all $n\geq n_0$ that
$$\theta_1\E^{\ell}_{\alpha+\e}(\mu^a_{n})+\theta_2\E^{\ell}_{\alpha-\e}(\mu^a_{n})\leq \E^{\ell}(\mu^a_{n}),$$
where $\E^\ell$, is the truncated energy.
Which gives
$$\theta_1\E^{\ell}_{\alpha+\e}(\mu^a_{n})+\theta_2\E^{\ell}_{\alpha-\e}(\mu^a_{n})\leq\frac{n(n-1)}{n^2}\tau_n
-\frac1{n^2}\sum_{i=1}^nV_n(\sigma^{(n, a)}_i)+\frac{\ell}{n},$$
Since $$V_n(x)\geq Q_{\alpha_n}(x),$$
and by a simple computation we can show that
$$\inf_{s\in\Bbb R_+}Q_{\alpha_n}(s)=\left\{\begin{aligned}&\alpha_n-\alpha_n\log\alpha_n,\;\;{\rm if}\;\alpha_n>0\;{\rm for\; some \;n}\\&0\qquad\qquad\qquad\;{\rm if}\;\alpha_n=0\;\forall\;n\end{aligned}\right.$$
Hence $$\theta_1\E^{\ell}_{\alpha+\e}(\mu^a_{n})+\theta_2\E^{\ell}_{\alpha-\e}(\mu^a_{n})\leq\frac{n(n-1)}{n^2}\tau_n
-\frac1{n}\alpha_n(1-\log\alpha_n)+\frac{\ell}n,$$
As $n$ goes to infinity, using the fact that $\alpha_n$ converge we obtain
$$\liminf_n(\theta_1\E^{\ell}_{\alpha+\e}(\mu^a_{n})+\theta_2\E^{\ell}_{\alpha-\e}(\mu^a_{n}))\leq\liminf\tau_n,$$
hence $$\theta_1\E^{\ell}_{\alpha+\e}(\sigma^a)+\theta_2\E^{\ell}_{\alpha-\e}(\sigma^a)\leq\liminf\tau_n,$$
applying the monotone convergence theorem, when $\ell$ goes to $+\infty$ we obtain
$$\theta_1\E_{\alpha+\e}(\sigma^a)+\theta_2\E_{\alpha-\e}(\sigma^a)\leq\liminf\tau_n,$$
Since $\theta_1\E_{\alpha+\e}(\sigma^a)+\theta_2\E_{\alpha-\e}(\sigma^a)=\E_\alpha(\sigma^a)$. It follows
$$\E_{\alpha}(\sigma^a)\leq\liminf\tau_n.$$
Furthermore $$\E^*_\alpha=\inf_{\nu\in\mathfrak{M}^1(\Bbb R)}\E_\alpha(\nu)\leq \E_\alpha(\sigma^a),$$
hence $$\E^*_\alpha\leq \E_\alpha(\sigma^a)\leq\liminf\tau_n.$$
Therefore $$\E^*_\alpha\leq  \E_\alpha(\sigma^a)\leq\liminf\tau_n\leq \limsup\tau_n\leq \E^*_\alpha,$$
in the last inequality we have used equation (\ref{eq5}).
Then we obtain $$\E_\alpha(\sigma^a)=\E^*_\alpha=\E_\alpha(\mu_\alpha^a).$$
This implies by unicity of the equilibrium measure that $\sigma^a=\mu_\alpha^a$. We have proved that $\mu_\alpha^a$ is
the only possible limit for a subsequence of the sequence $(\mu_n^a)$. It follows
that the measure $\mu_n^a$ itself converges:
for all bounded continuous function on $[a,+\infty[$, $$\lim_{n\to\infty}\int_a^{+\infty}f(x)\mu_n^a(dx)=\int_a^{+\infty}f(x)\mu_\alpha^a(dx),$$
and $$\lim_{n\to\infty}\tau_n=\E^*_\alpha.$$
\qquad\qquad\qquad\qquad\qquad\qquad\qquad\qquad\qquad\qquad\qquad\qquad\qquad\qquad\qquad\qquad\qquad$\blacksquare$
\begin{proposition}---
  For $a\geq 0$, let $$E^*_n(a)= E\Big(x_1^{(n, a)},...,x_n^{(n, a)}\Big),$$ be the minimum of the energy at the zeros of $H_n^{\lambda_n, a}$. If $\displaystyle\lim\limits_{n\to+\infty}\frac{\lambda_n}{n}=0$, then for $n$ large enough
 $$E_n^*(a)=\E^*n^2-n\Big(\lambda_n+\frac n2\Big)\log n+o(n^2),$$
 where $\E^*$ is the energy given in proposition 3.5.

\end{proposition}
For $a=0$, one obtains $$E_n^*=\Big(\frac34+\frac12\log6\Big)n^2-n\Big(\lambda_n+\frac n2\Big)\log n+o(n^2).$$
For $\lambda_n=0$, the last asymptotic formula gives an approximation of the energy in the case of half range Hermite polynomials.
$$E_n^*\sim-\frac {n^2}2\log n.$$
Such result is proved in \cite{ch}.

More general for $\lim\limits_{n\to+\infty}\frac{\lambda_n}{n}=\alpha>0$, the energy $\E^*_\alpha$, has been computed see for instance \cite{bo}. In this case one can easily obtains as $n\to+\infty$,
$$E_n^*=\E^*_\alpha n^2-n\Big(\lambda_n+\frac n2\Big)\log n+o(n^2).$$
Moreover, it has been proved in \cite{bo}, that for $a=0$, and $\alpha$ small enough
$$\lim_{n\to+\infty}\Big(\frac1{n^2}E_n^*+\Big(\frac{\lambda_n}{n}+\frac12\Big)\log n\Big)=\Big(\frac34+\frac12\log 6+C\alpha\Big)+o(\alpha),$$
where  $\displaystyle C=\frac1{432}\Big(-36 (-6 + \sqrt6) + (54 - 161 \sqrt6) \log2 +
    27 (10 + \sqrt6) \log3\Big)\approx0.6045$, and $o(\alpha)$ is a small term in $\alpha$.

{\bf Proof.}---
We saw $$E_n^*(a)=E(\sqrt n\sigma_1^{(n, a)},...,\sqrt n\sigma^{(n, a)}_n),$$
hence $$E^*_n(a)=\sum_{i<j}\log\frac{1}{\Big|\sqrt n\sigma_i^{(n, a)}-\sqrt n\sigma_n^{(n, a)}\Big|}+\sum_{i=1}^nV(\sqrt n\sigma_i^{(n, a)})-n\log(C_{\lambda_n}),$$
since $$V(\sqrt n x)=nx^2+2\lambda_n\log\frac1{\sqrt nx}+\log\Big(2+\frac{2b_n}{\sqrt n x}+\frac{\beta_n}{\sqrt n x-a}\Big)+\log(C_{\lambda_n}),$$
hence $$V(\sqrt n x)=nV_n(x)-\lambda_n\log n+\log(C_{\lambda_n}),$$
where $\displaystyle V_n(x)=\Big(x^2+2\alpha_n\log\frac1x+\frac1n\log\big(2+\frac{2b_n}{\sqrt n x}+\frac{\beta_n}{\sqrt n x-a}\big)\Big).$
Hence $$E^*_n(a)+\Big(n\lambda_n+\frac{n(n-1)}{2}\Big)\log n=E_{\alpha_n}^*,$$
By the previous theorem we get
$$\lim_{n\to+\infty}\frac1{n^2}\Big(E^*_n(a)+\Big(n\lambda_n+\frac{n(n-1)}{2}\Big)\log n\Big)=\lim_{n\to+\infty}E_{\alpha_n}^*=\E_\alpha,$$
Hence $$\lim_{n\to+\infty}\Big(\frac1{n^2}E^*_n(a)+\Big(\frac{\lambda_n}{n}+\frac{1}{2}\Big)\log n\Big)=\E_\alpha,$$
and the case $\alpha=0$, gives the desired result.\qquad\qquad\qquad\qquad\qquad\qquad\qquad\qquad$\blacksquare$
\begin{center}{\bf Plotting of density of zeros.}\end{center}

\begin{figure}[h]
\centering\scalebox{0.5}{\includegraphics[width=17cm, height=10cm]{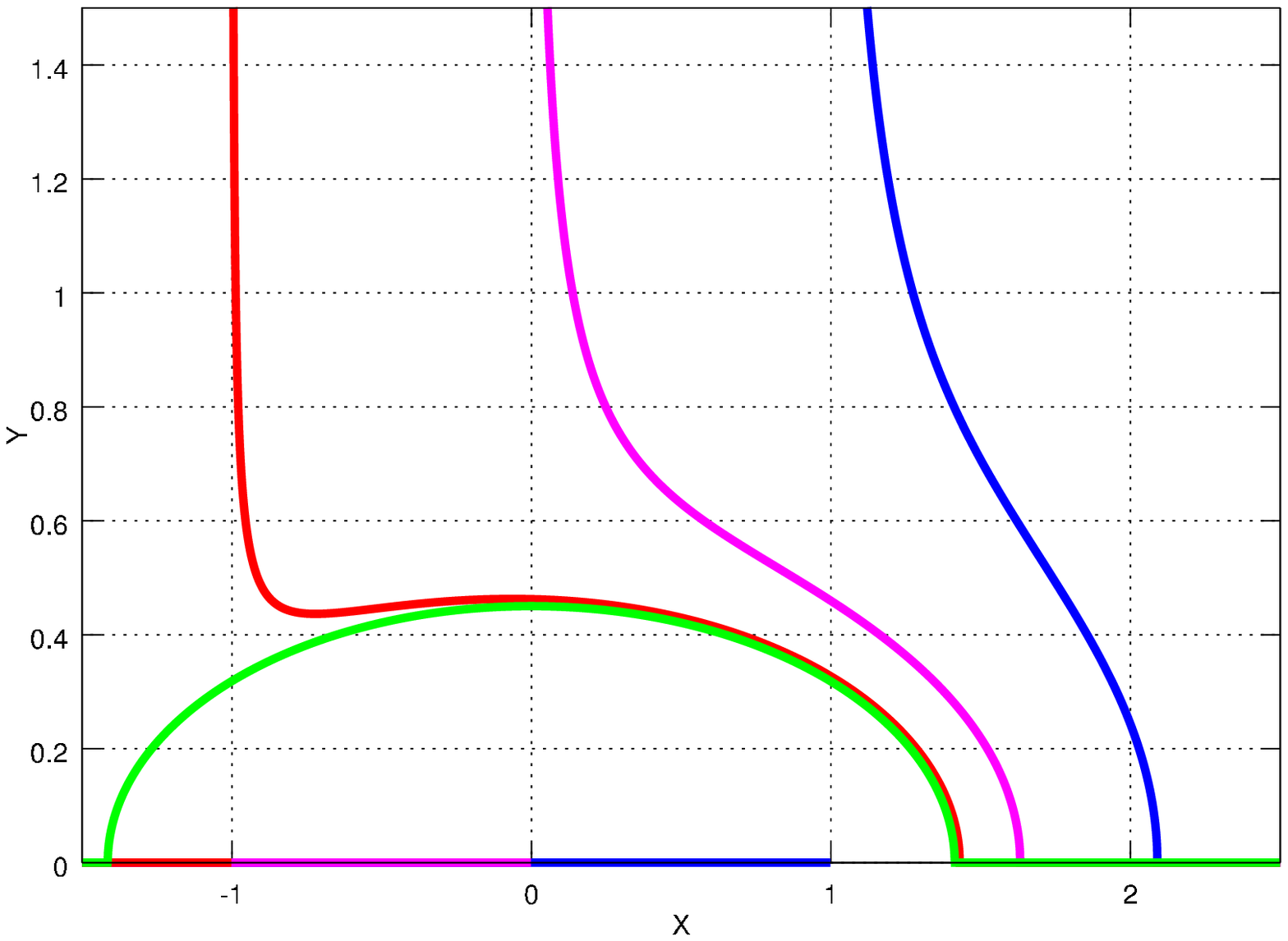}}
\caption{\bf Plot of $f_{\alpha, a}$ for $\alpha=0$}
 \end{figure}
 \newpage
 \begin{flushleft}
\textcolor{blue}{\rule{1cm}{1pt}} for  $a=1$, $b\approx 2.09$, $\sigma=1$, density of eigenvalues in $\Sigma_\sigma=[1,+\infty[$.\\
\textcolor{OliveGreen}{\rule{1cm}{1pt}} for  $a=0$, $b=\frac23\sqrt6\approx 1.632$, $\sigma=0$, density of zeros in $\Sigma_\sigma=[0,+\infty[$.\\
\textcolor{red}{\rule{1cm}{1pt}}for  $a=-1$, $b=1.43$, $\sigma=-1$, density of zeros in $\Sigma_\sigma=[-1,+\infty[$.\\
\textcolor{BlueGreen}{\rule{1cm}{1pt}}for  $a=-\sqrt 2$, $b=\sqrt 2$, $\sigma\leq-\sqrt 2$, density of zeros in $\Sigma_\sigma=[\sigma,+\infty[$.\\
\end{flushleft}

\begin{figure}[h]
\centering\scalebox{0.6}{\includegraphics[width=15cm, height=10cm]{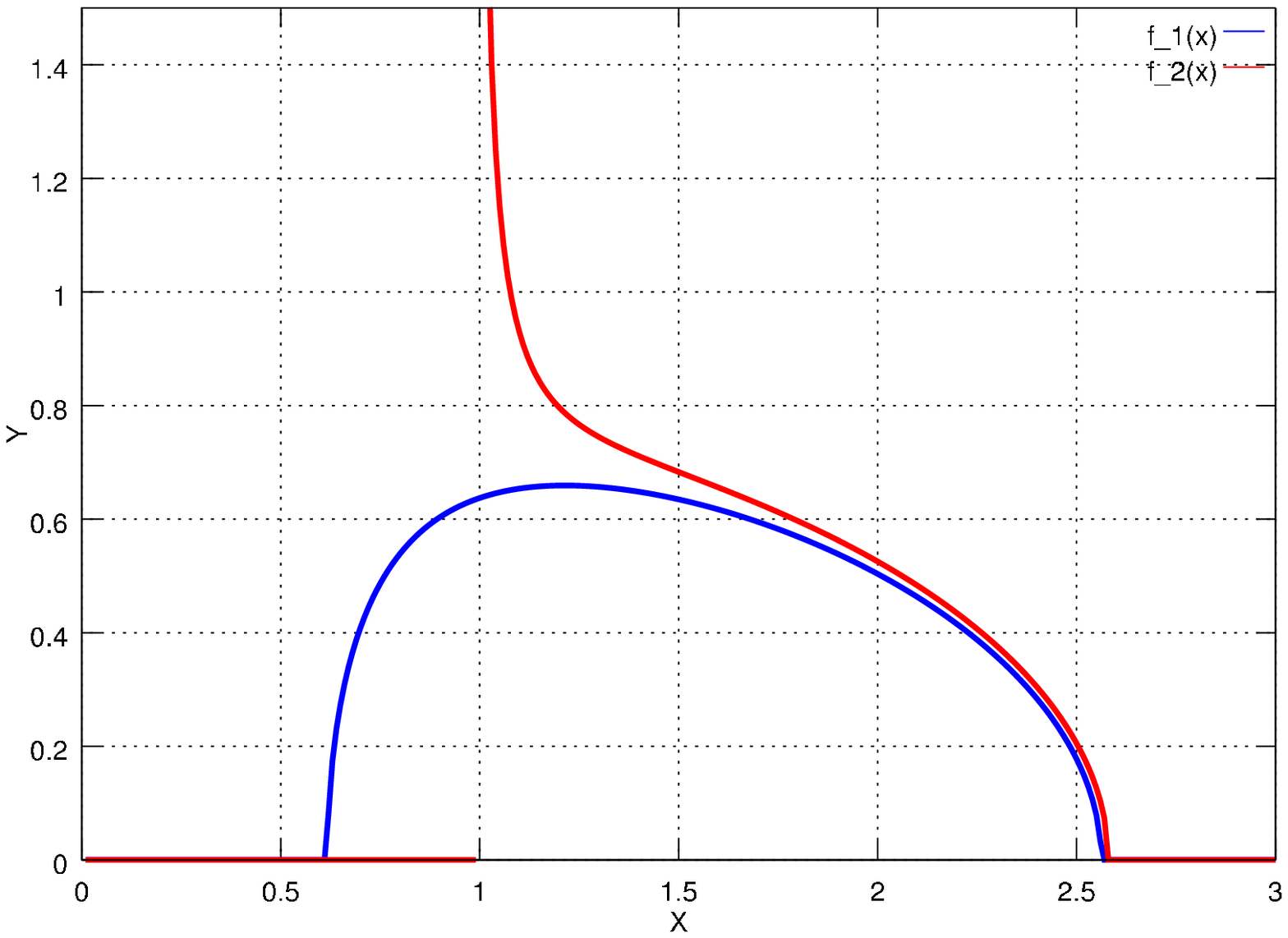}}
\caption{\bf Plot of $f_{\alpha, a}$ for $\alpha=2$.}

\end{figure}
 \begin{flushleft}
 \textcolor{red}{\rule{1cm}{1pt}} $a=1$, $b=2.58$, density of zeros in $[1,+\infty[$.\\
 \textcolor{blue}{\rule{1cm}{1pt}} $a=a_c\approx 0.618$, $b=b_c\approx 2.562$,
density of zeros in $]a,+\infty[$, for $0\leq a\leq a_c$
 \end{flushleft}

\section{ Truncated Generalized Hermite polynomials}
Let $a\in\Bbb R$, and $\lambda\geq 0$. Consider in $]-\infty,-a]\cup[a,+\infty[$, the probability measure
$$w_\lambda(x)=\frac{1}{C_{\lambda, a}}|x|^{2\lambda}e^{-x^2}=e^{-v_\lambda(x)},$$
where $\displaystyle C_{\lambda, a}=2\int_a^{+\infty}|x|^{2\lambda}e^{-x^2}dx$. For $a=0$, $C_\lambda=\Gamma(\lambda+\frac12)$. As in the first section to the weight function $w_\lambda$ we associate the sequence of orthonomalizing polynomials $H^{\lambda, a}_n$.
 They satisfies $$\int_{]-\infty,-a]\cup[a,+\infty[}H^{\lambda, a}_n(x)H_m^{\lambda, a}(x)w_\lambda(x)dx=\delta_{nm}.$$
 By the symmetry of the weight they satisfies the three-terms recurrence relation
$\displaystyle H_0^{\lambda, a}(x)=1$, $\displaystyle H_1^{\lambda, a}(x)=\frac{x-a_1}{b_1}$ and for all $n\geq 1$,
$$xH^{\lambda, a}_n(x)=a_{n+1}H^{\lambda, a}_{n+1}(x)+a_nH^{\lambda, a}_{n-1}(x).$$
Moreover if we write $H^{\lambda, a}_n(x)=\gamma_nx^n+$ lower order terms, then $\displaystyle a_n=\frac{\gamma_{n-1}}{\gamma_n}$.

The sequences $a_n, b_n$  and $\gamma_n$ all depend in the parameters $\lambda$ and $a$, but we omit this dependance in the notations.
\begin{proposition}---
Let $b_n=\max(|a_n|, |a_{n+1}|)$, then
$$a_n^2+a_{n+1}^2=n+\lambda+\frac12+a(H_n^{\lambda, a}(a))^2w_\lambda(a),$$
and for $a\neq 0$, $$(H_n^{\lambda, a}(a))^2w_\lambda(a)=O(b_n^2).$$
Moreover $$\lim_{n\to +\infty}b_n^2=+\infty.$$
If an addition $a=0$ then $a_n=O(\sqrt{n+\lambda})$.
\end{proposition}
{\bf Proof.}---By orthogonality we have
$$n=\int_{]-\infty, -a]\cup[a,+\infty[}xH_n^{\lambda, a}(x)'H_n^{\lambda, a}(x)w_\lambda(x)dx,$$
Integrate by part we obtain
$$n=-a(H_n^{\lambda, a}(a))^2w_\lambda(a)-\frac12-\frac12\int_{]-\infty, -a]\cup[a,+\infty[}x\big(H_n^{\lambda, a}(x)\big)^2w'_\lambda(x)dx,$$
Since $$w'_\lambda(x)=\Big(-2x+\frac{2\lambda}{x}\Big)w_\lambda(x),$$
it follows that $$\int_{]-\infty, -a]\cup[a,+\infty[}x\big(H_n^{\lambda, a}(x)\big)^2w'_\lambda(x)dx=-2\big(xH_n^{\lambda, a}\mid xH_n^{\lambda, a}\big)+2\lambda\big(H_n^{\lambda, n}\mid H_n^{\lambda, a}\big),$$
From the recurrence relation one gets,
$$n=-a(H_n^{\lambda, a}(a))^2w_\lambda(a)-\frac12+a_n^2+a_{n+1}^2-\lambda,$$
and the result follow.
For the second step one can see that
$$\frac{(H_n^{\lambda, a}(a))^2w_\lambda(a)}{b_n^2}\leq \frac{1}{a},$$
moreover, $b_n^2\geq n+1$, hence $\lim_{n\to+\infty}b_n=+\infty$.\\
If $a=0$, then
$$a_n^2+a_{n+1}^2=2n+1+2\lambda,$$
and $$|a_n|\leq\sqrt{2n+2\lambda+1}.$$
Which gives the desired result.
\begin{proposition}---
The polynomials satisfies the differential equations
\begin{description}
\item [(1)]
$$H^{\lambda, a}_n(x)'=A_n(x)H^{\lambda, a}_{n-1}(x)-B_n(x)H^{\lambda, a}_n(x),$$
where $$A_n(x)=2a_n\Big(1+\frac{aH^{\lambda, a}_n(a)^2w_\lambda(a)}{x^2-a^2}\Big),$$
and $$B_n(x)=2a_nH^{\lambda, a}_n(a)H^{\lambda, a}_{n-1}(a)\frac{x}{a^2-x^2}+\frac{2\lambda a_n}{x}\int_{|y|\geq a}
H^{\lambda, a}_n(y)H^{\lambda, a}_{n-1}(y)w_\lambda(y)\frac{dy}{y}.$$
\item[(2)]
$$H^{\lambda, a}_n(x)''+R_n(x)H^{\lambda, a}_n(x)'+S_n(x)H^{\lambda, a}_n(x)=0,$$
where $$R_n(x)=-2x+\frac{2\lambda}{x}-\frac{A_n'(x)}{A_n(x)},$$ and
$$S_n(x)=B'_n(x)-B_n(x)\frac{A'_n(x)}{A_n(x)}-B_n(x)\Big(2x-\frac{2\lambda}{x}+B_n(x)\Big)+\frac{a_n}{a_{n-1}}A_n(x)A_{n-1}(x)$$
\end{description}
\end{proposition}
\section{Density of zeros of trauncated symmetric generalized Hermite polynomials}
Let $$V(x)=x^2+2\lambda\log\frac1{|x|}+\log C_\lambda+\log\Big(2+2w_\lambda(a)(H^{\lambda, a}_n(a))^2\frac{a}{x^2-a^2}\Big)$$
Consider the system
of $n$ movable unit charges in $\Big(]-\infty, -a[\cup]a; +\infty[\Big)^n$ in the presence of the external potential $V(x)$. If $x_1,...,x_n$ are the position of the particles arranged in decreasing order. The total energy of the system is
 $$ F(x_1,\cdots,x_n)=2\sum_{i< j}\log\frac{1}{|x_i-x_j|}+\sum_{i=1}^nV(x_i),$$
 \begin{proposition}---The energy $F(x_1,...,x_n)$ is minimal at the zeros of the truncated generalized Hermite polynomials $H^{\lambda, a}_n$.
 \end{proposition}
 See for instance \cite{I}.

 Let $x^{(n, a)}_1,...,x^{(n, a)}_n\in]-\infty, -a]\cup [a, +\infty[$ denote the zeros of truncated generalized Hermite polynomials $H^{\lambda_n, a}_n$, where $\lambda_n$ be some positive real sequence, and defined on $]-\infty, -a]\cup[a, +\infty[$  the probability measure
$$\nu_n^a=\frac1n\sum_{i=1}^n\delta_{\sigma^{(n, a)}_i}, $$
where $\displaystyle \sigma^{(n, a)}_i=\frac{x^{(n, a)}_i}{b_n}$ are the rescaled zeros. Moreover consider the modified rescaled energy in the sense
 $$E(x_1,\cdots,x_n)=F(x_1,...,x_n)-n\log C_{\lambda_n}=2\sum_{i< j}\log\frac{1}{|x_i-x_j|}+b_n\sum_{i=1}^nV_n(x_i),$$
 where $\displaystyle V_n(x)=x^2+2\alpha_n\log\frac1x+\frac1{b_n}\log\Big(2+\frac{2a\beta_n}{b_n^2 x^2-a^2}\Big)$, $\displaystyle\alpha_n=\frac{\lambda_n}{b_n}$, and $\displaystyle\beta_n=w_{\lambda_n}(a)(H_n^{\lambda_n}(a))^2$. Then the minimum of the energy is attained at the rescaled zeros $\displaystyle \sigma_i^{(n, a)}$ of the polynomials $H_n^{\lambda_n, a}$.
 For $a=0$ one can take instead of $b_n$ the sequence $n$.
\begin{theorem}---
We assume that, the sequence $\alpha_n$ converge to some $\alpha\geq 0$. Then for all $a\geq 0$, the measure $\nu_n^a$ converge for the tight topology to some probability measure $\nu_\alpha^\sigma $ supported by $[-b, -\sigma]\cup[\sigma,b]$ with density
$f_{\alpha, \sigma}$, such that
 \begin{description}
\item[(1)] If $a=0$ and $\alpha>0$, then there is a unique $\sigma_0=\sigma(\alpha)>0$ and  a unique $b_0=b(\alpha)$ such that $b>\sigma_0>0$, and
$$f_{\alpha, \sigma_0}(x)=\frac1{\pi|x|}\sqrt{(b_0^2-x^2)(x^2-\sigma_0^2)},$$
where $$\displaystyle b_0=\sqrt{1+\alpha+\sqrt{1+2\alpha}},$$ and
$$\sigma_0=\sqrt{1+\alpha-\sqrt{1+2\alpha}}.$$

\item[(2)] If $0<a\leq\sigma_0$ and $\alpha>0$, then $\sigma(a, \alpha)=\sigma_0$ and the density on $\Sigma$ still the same as in {\bf (1)}.\\
 \item[(3)]If $a>\sigma_0$ and $\alpha>0$, then $\sigma(a, \alpha)=a$ and,
  $$f_{\alpha, \sigma}(x)=\frac1{\pi|x|}\sqrt{\frac{b^2-x^2}{x^2-\sigma^2}}\Big(x^2-\frac{\alpha \sigma}{b}\Big).$$
  where  $b=b(\alpha, \sigma)$ is the unique solutions of the following equations
$$\sigma b-\alpha\geq 0,\;and\; \frac{b^2}{2}+\frac{\alpha\sigma}{b}-\frac{\sigma^2}{2}-\alpha-1=0.$$
 \item[(4)]If $a\geq0$ and $\alpha=0$, then $\sigma(a, 0)=a$ and,
     $$f_{0, \sigma}(x)=\frac{|x|}{\pi}\sqrt{\frac{2+\sigma^2-x^2}{x^2-\sigma^2}},$$
     where $b=\sqrt{\sigma^2+2}$
    \end{description}
\end{theorem}
One recover's from {\bf (1)} or {\bf (4)} the density of zeros of classical Hermite polynomials which is given by the semi-circle law $$f_{0,0}(x)=\frac{1}{\pi}\sqrt{2-x^2}.$$
Moreover using mathematica one can find the explicit value of $b$ in {\bf (4)}, but it is very complicate to put there here.

The proof of theorem 4.2 follow as theorem 2.2 with slight modifications.

As in the previous section, in \cite{bo} the energy $E^*_{\alpha, a}$ of the measure $\nu_\alpha^a$ in the previous theorem has been computed, and as in the previous section one can prove that, the asymptotic of the equilibrium energy of the zeros is as follows:
$$E_n^*(a)=E^*_{\alpha, a}n^2-n\Big(\lambda_n+\frac n2\Big)\log n+o(n^2).$$
For the particular case, if $\displaystyle \lim_{n\to +\infty}\frac{\lambda_n}{n}=\alpha=0$, one gets
$$E_n^*(a)=\Big(\frac34+\frac12\log 2+a^2\Big)n^2-n\Big(\lambda_n+\frac n2\Big)\log n+o(n^2),$$

\begin{center}{\bf Plotting of the density of zeros $f_{\alpha, \sigma}$}\end{center}
\begin{figure}[h]
\centering\scalebox{0.6}{\includegraphics[width=17cm, height=8cm]{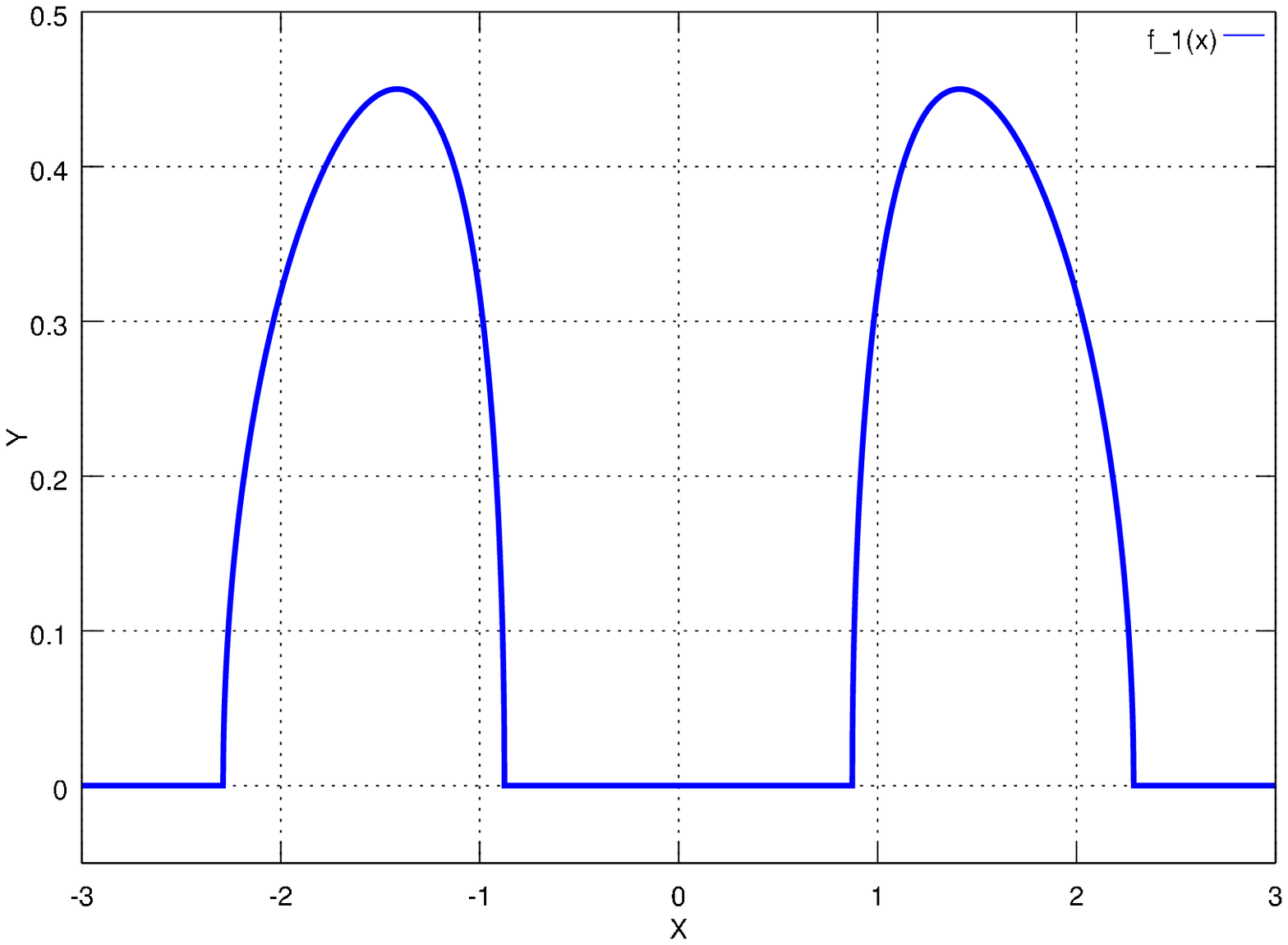}}
\caption{\bf $\alpha=2$, $\sigma=a_c\approx 0.874$, $b=b_c\approx 2.2882$}
Density of zeros on $]-\infty,-a]\cup[a,+\infty[$, for all $0<a\leq 0.874$.
\end{figure}
\newpage

\begin{figure}[h]
\centering\scalebox{0.5}{\includegraphics[width=17cm, height=10cm]{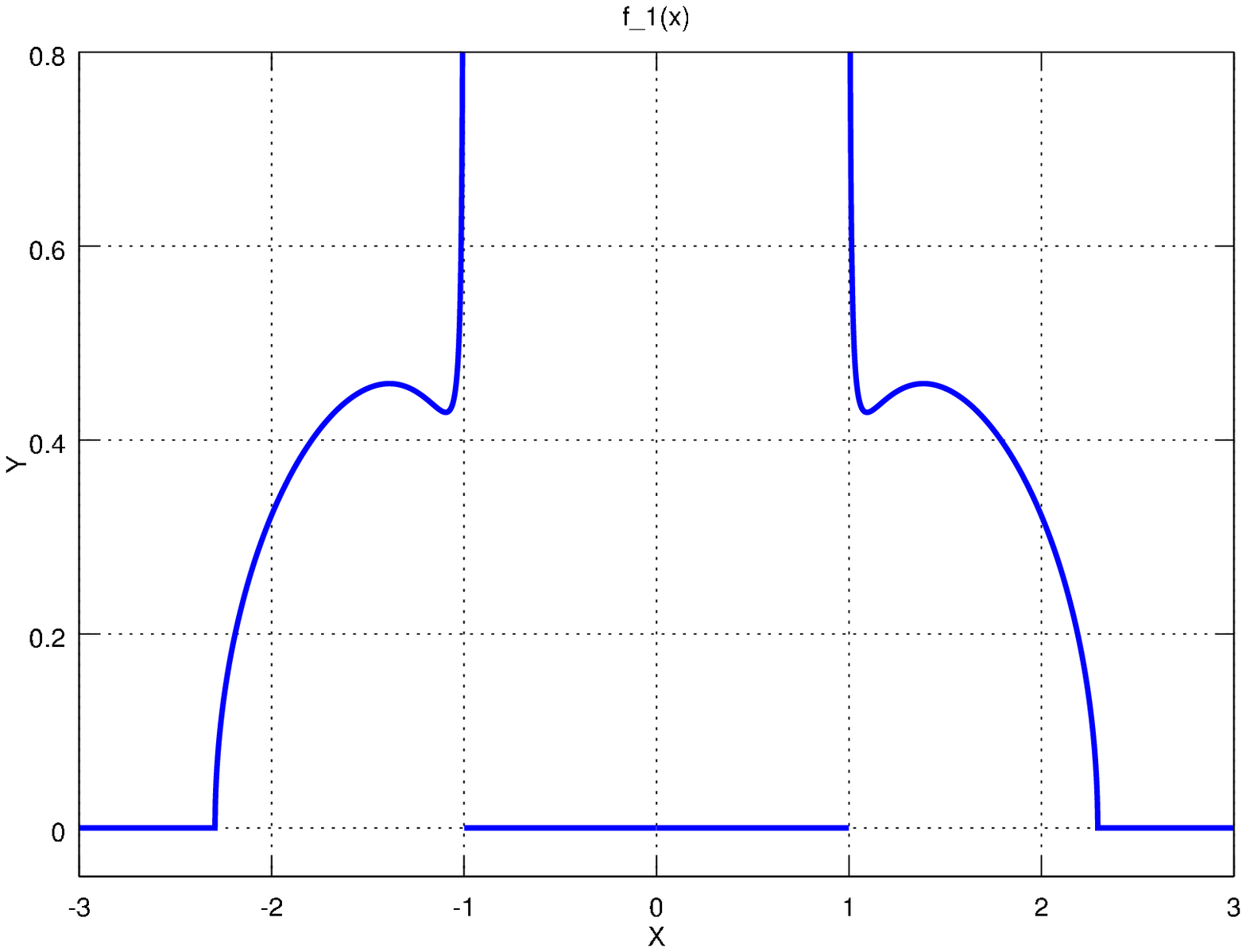}}
\caption{\bf $\alpha=2$, $a=1>a_c$, $b=2.2924$}
Density of zeros on $]-\infty,-1]\cup[1,+\infty[$.
\end{figure}

\begin{figure}[h]
\centering\scalebox{0.6}{\includegraphics[width=15cm, height=10cm]{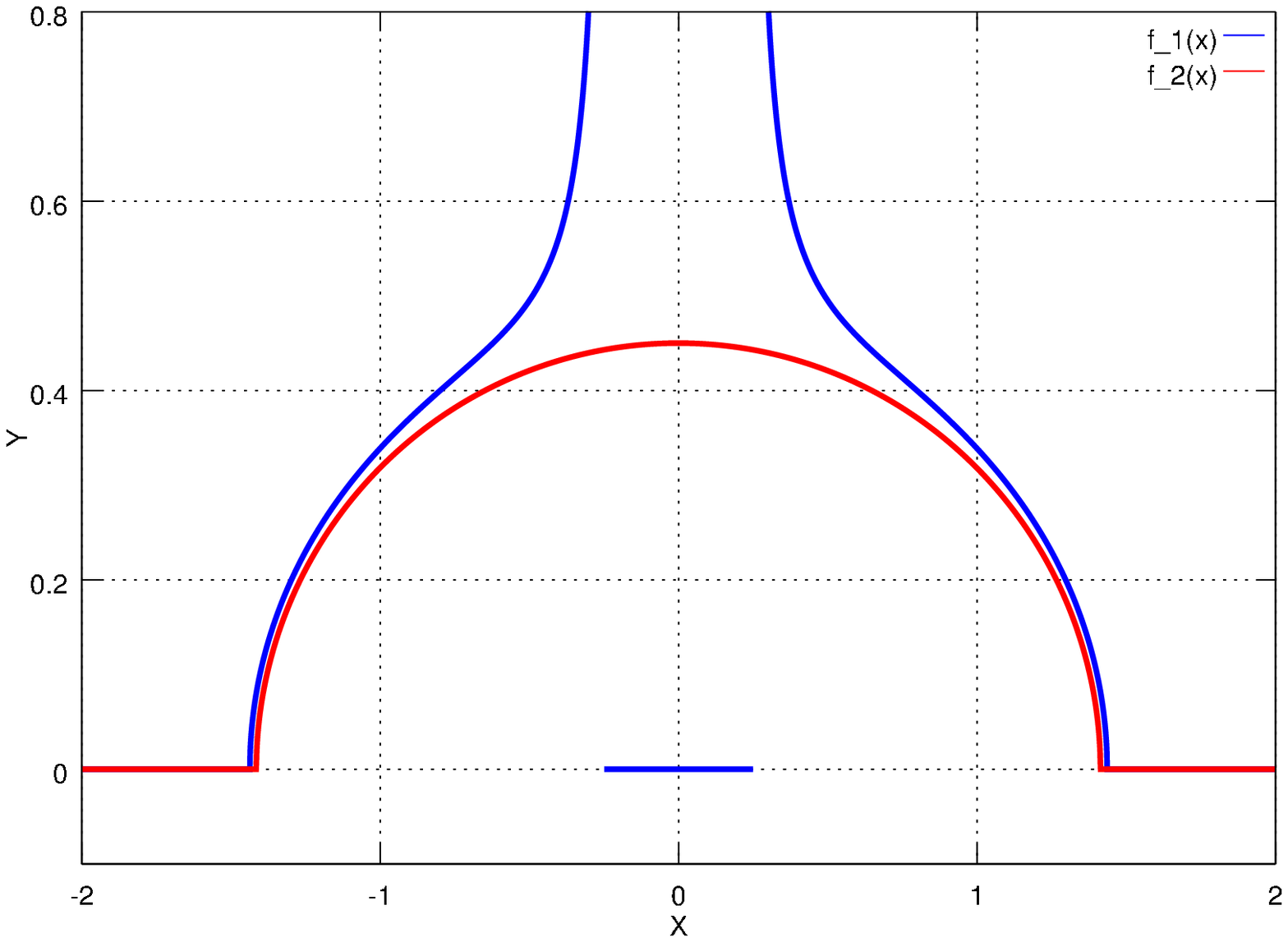}}
\caption{\bf $\alpha=0$}
\begin{flushleft}\textcolor{red}{\rule{1cm}{1pt}}  $a=a_c=0$, $b=b_c=\sqrt 2$ semi-circle law, density of zeros on $\Bbb R$.\\
\textcolor{blue}{\rule{1cm}{1pt}} $a=0.25$, $b= 1.436$ density of zeros on $]-\infty,-0.25]\cup[0.25,+\infty[$.
\end{flushleft}

\end{figure}

Address:  College of Applied Sciences
   Umm Al-Qura University
  P.O Box  (715), Makkah,
  Saudi Arabia.\\
E-mail: bouali25@laposte.net \& mabouali@uqu.edu.sa


\newpage


\begin{thebibliography}{9}

\bibitem{B} W. Bauldry, \emph{Estimates of asymmetric Freud polynomials on the real line}, J. Approximation
Theory, 63, 225-237, (1990).

\bibitem{bo} M. Bouali, \emph{Density of Positive Eigenvalues of the Generalized Gaussian Unitary Ensemble.} arXiv:1409.0103v1 [math.PR], (submitted).
\bibitem{boo}--------, \emph{Generalized $\beta$-Gaussian Ensemble Equilibrium measure method.} arXiv:1409.0126v1 [math.PR], (Submitted).
\bibitem{C} S.S. Bonan \& D.S. Clark, \emph{Estimates of the Hermite and the Freud polynomials}, J.
Approximation Theory, 63, 210-224, (1990).
\bibitem{IC} Y. Chen \& M.E.H. Ismail, \emph{Ladder operators and differential equations for orthogonal
polynomials}, J. Phys. A30, 7818-7829, (1997).
\bibitem{ch} Y. Chen, M.E.H. Ismail and W. Van Assche, \emph{Tau-function constructions of the recurrence
coeffcients of orthogonal polynomials}, Advances in Appl. Math., 20 (1998),
141-168
\bibitem{F} J. Faraut. \emph{ Logarithmic potential theory, orthogonal polynomials,
and random matrices}, CIMPA Scool, (2011).
 \bibitem{SD}  D. S. Dean\& S. N. Majumdar,\emph{ Extreme Value Statistics of Eigenvalues of Gaussian Random Matrices.} Physical Review E: Statistical, Nonlinear, and Soft Matter Physics 77, 4, (2008).
\bibitem{S}  S. N. Majumdar, C. Nadal, A. Scardicchio. \& P. Vivo. \emph{ How many eigenvalues of a Gaussian random matrix are positive?} Physical Review E 83, 041105 (2011).

 \bibitem{I} M. Ismail, \emph{An Eectrostatics Model For Zeros Of General Othogonal Polynomials}. Pacific Journal of Mathematics, Vol. 193, No. 2, (2000).
     \bibitem{sz} G. {Szeg}'o, \emph{Orthogonal Polynomials}, Fourth Edition, Amer. Math. Soc., Providence,
(1975).
     \bibitem{st} T. J. Stieltjes, \emph{Sur quelques th\'eor$\grave{e}$mes d'alg$\grave{e}$bre}, Comptes Rendus de l'Academie des
Sciences, Paris, 100, 439-440; Oeuvres Compl$\grave{e}$tes, Vol. 1, 440-441, (1885).
\bibitem{stl}--------,Sur les polyn\^omes de Jacobi, Comptes Rendus de l'Academie des Sciences,
Paris, 100, 620-622; Oeuvres Compl$\grave{e}$tes, Vol. 1, 442-444, (1885).

\end{thebibliography}
\end{document}